\documentclass[a4paper,10pt]{article}
\usepackage{graphicx}
\usepackage{amsmath}
\usepackage{amssymb}
\usepackage{hyperref}
\usepackage{enumitem}
\usepackage{mathtools}

\usepackage{amsthm}

\usepackage{siunitx} % For units
\usepackage[numbers]{natbib}

\usepackage{multicol} 
\usepackage{float}

\usepackage{titlesec}

\titleformat{\section}{\normalfont\Large\bfseries\MakeUppercase}{\thesection}{1em}{}
\titleformat{\subsection}{\normalfont\large\bfseries\MakeUppercase}{\thesubsection}{1em}{}
\titleformat{\subsubsection}{\normalfont\normalsize\bfseries\MakeUppercase}{\thesubsubsection}{1em}{}

\usepackage[top=1in, bottom=1in, left=0.75in, right=0.75in]{geometry}

\usepackage{helvet}

\usepackage{fancyhdr}
\makeatletter
\patchcmd{\@maketitle}{\thispagestyle{plain}}{\thispagestyle{fancy}}{}{}
\makeatother

\fancypagestyle{plain}{ % redefine plain so first page also uses same style
  \fancyhf{}
  \fancyhead[L]{Proceedings of the URECA@NTU 2024-25}
  \fancyfoot[C]{\thepage}
}

\makeatletter
\renewenvironment{abstract}{
    \small
    \noindent\textbf{\textit{\abstractname}}.\hspace{1em}} % Bold and italic "Abstract" followed by a dot and space
    {}
\makeatother
\usepackage[parfill]{parskip} % Removes paragraph indent and adds spacing between paragraphs
\usepackage{setspace}

\title{A Brief Review of Fixed Points, Hex Game and Hex Theorem}
\date{}
\author{
    {Russell Cho Yang} \\
    School of Physical and Mathematical Sciences,
    Nanyang Technological University, Singapore
}

\begin{document}

%%%% Uncomment and insert the date if needed
% \date{}
\maketitle
\begin{multicols}{2}

\begin{abstract}
The Kakutani Fixed Point Theorem is a cornerstone of Game Theory, showing how equilibrium points emerge in $n$-person games. This research looks at the theorem through the rules and strategies of the board game Hex. Inspired by the impossibility of ties in certain games, we explore how the game of Hex is equivalent to Brouwer Fixed Point Theorem. Kakutani’s Fixed Point Theorem is a generalization of Brouwer’s Fixed Point Theorem and Sperner's lemma, on the other hand, is a combinatorial result concerning the labeling of the vertices of simplices and their triangulations. Sperner's lemma is later known to be equivalent to Brouwer's fixed point Theorem. Due to this, By analyzing the game of Hex finite moves and stable states (Nash equilibria), we uncover how Hex can illustrate deep concepts like convexity and continuity from its connection with Sperner's lemma. This report will focus more on Hex and its equivalence to the Brouwer fixed point Theorem. The detailed analysis of Sperner's lemma as well as the proof showing it is equivalent to Kakutani Fixed Point Theorem or whether Sperner's Lemma is equivalent to Brouwer Theorem, will not be covered in this report, but some of its key definitions and results can be found in the Appendix.
\end{abstract}

% Include at least three keywords of phrases
\textbf{Keywords:} Hex Boardgame, Brouwer Fixed Point, Kakutani Fixed Point

\section{Introduction}
Hex is a two-player abstract strategy game played on a diamond-shaped board made up of hexagonal cells, usually arranged in an 11-by-11 grid. The two players, commonly referred to as Red and Blue, take turns marking empty hexes with their color. The objective is to create an unbroken path of one's own color connecting their assigned sides of the board—Red aims to connect the top and bottom, while Blue tries to connect the left and right edges. An example of a winning setup for Blue is shown in Figure 1 \cite{Liu}.

The Hex Theorem states that a game of Hex can never end in a draw—one player must always win. This is because the only way to block the opponent from completing their path is to complete your own first. Although the result seems obvious, formally proving it can involve deep concepts from topology. Mathematician David Gale provided a simpler proof using graph theory, and also demonstrated that the Hex Theorem is logically equivalent to the Brouwer Fixed Point.

\section{Hex Board Game}

\begin{figure} [H]
    \centering
    \includegraphics[width=\linewidth]{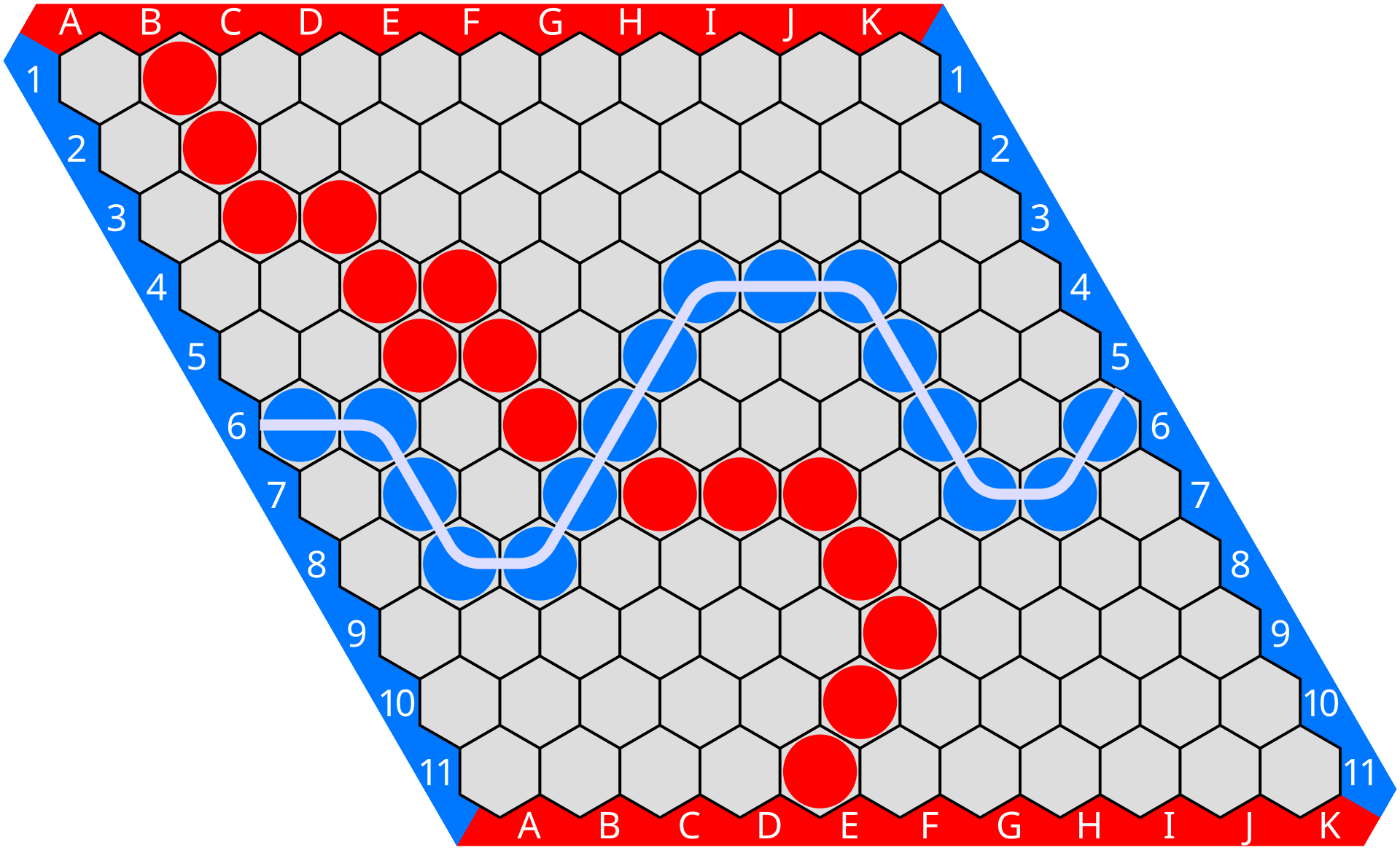}
    \caption{$11 \times 11$ Hex Configuration Win For Blue}
    \label{fig:hex_blue_wins}
\end{figure}

For $11\times11$ Hex board as shown in Figure 1 \cite{Liu}, it is ultra-weakly solved. This means the first player always wins in perfect-play and this can be proved using a simple strategy-stealing argument.
First, we state two lemmas:

\textbf{Lemma 1}: Having extra/random pieces of your own color lying on the
board cannot hurt you.

\textbf{Lemma 2}: Hex cannot end in a draw.

To prove the first lemma, let's say an extra piece at position x on the board. If this extra piece's position is part of your winning strategy, then during your turn when you should be playing at position x, you could instead lay down another piece somewhere else. If x is not part of your winning strategy, then it would not matter if it is occupied.

For the second lemma, we will first take a look at an intuitive proof, while the rigorous proof will be shown in the next section. For the intuitive proof, first, take a board and put Red on the two sides of Red and Blue on the two sides of Blue. Put the rest of the stones on the board somewhere.

Between any two cells that touch with different colors, draw a line segment between them as seen in Figure 2 \cite{Hayward}. Observation about line segments:

- We will always have Blue on the left of the line and Red on the right

- Line segment path cannot end in the middle of the board

- Line segment path must end in one of the four corners (somewhere on the outside of the board)

There are two places for the line segment path to end up (starting from the Top-Left Corner):

- Surrounding Blue (Bottom Left); Red wins, or

- Surrounding Red (Top Right) - Blue wins

This shows only one player will have a full path of same colored connected hexes from their side to their opposite side. This is due to the fact that the two connection directions (Left–Right and Top–Bottom) are orthogonal, they cannot both happen without intersecting — which Hex forbids (no shared hexes).

\begin{figure} [H]
    \centering
    \includegraphics[width=\linewidth]{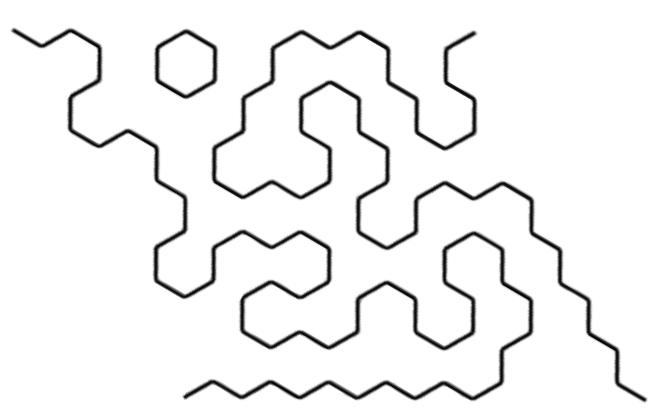}
    \caption{Hex Board Line Segments}
    \label{fig:no_draws}
\end{figure}

To prove that the first player always has a winning strategy in Hex, imagine two players, A (who goes first) and B (who goes second). Suppose, for sake of contradiction, that B has a guaranteed way to win. Then A could start the game by placing a random move anywhere, and from that point onward, follow B’s winning strategy — treating themselves as if they were the second player. According to a known result (Lemma 1), having that one extra piece from the random opening move cannot hurt A's chances. So A would be able to use B's strategy to win — meaning B couldn’t have had a winning strategy after all.

Also, since Hex cannot end in a draw (Lemma 2), and B can't win, it must be A who wins. Therefore, the first player must always have a winning strategy in Hex.

\section{Brouwer Fixed Point Theorem}
From Figure 3 \cite{Richard}, we can prove Brouwer's Fixed Point Theorem in the 1-Dimensional case. The continuous function $f$ is defined on a closed interval $[a, b]$ and takes values in the same interval. The statement that this function has a fixed point amounts to saying that its graph (dark green) intersects that of the function defined on the same interval $[a, b]$, which maps $x$ to $x$ (light green). Intuitively, any continuous line from the left edge of the square to the right edge must necessarily intersect the green diagonal. To prove this, consider the function $g$ which maps $x$ to $f(x) - x$. Observe that $g(0) = f(0) - 0 = f(0) \geq 0$ since $f(0) \in [0,1]$. At the same time, $g(1) = f(1) - 1 \leq 0$. Therefore, as $g(x)$ is $\geq 0$ on $a=0$ and $\leq 0$ on $b=1$. By the intermediate value theorem, $g$ has a zero in $[a, b]$; this zero is a fixed point.

\begin{figure} [H]
    \centering
    \includegraphics[width=\linewidth]{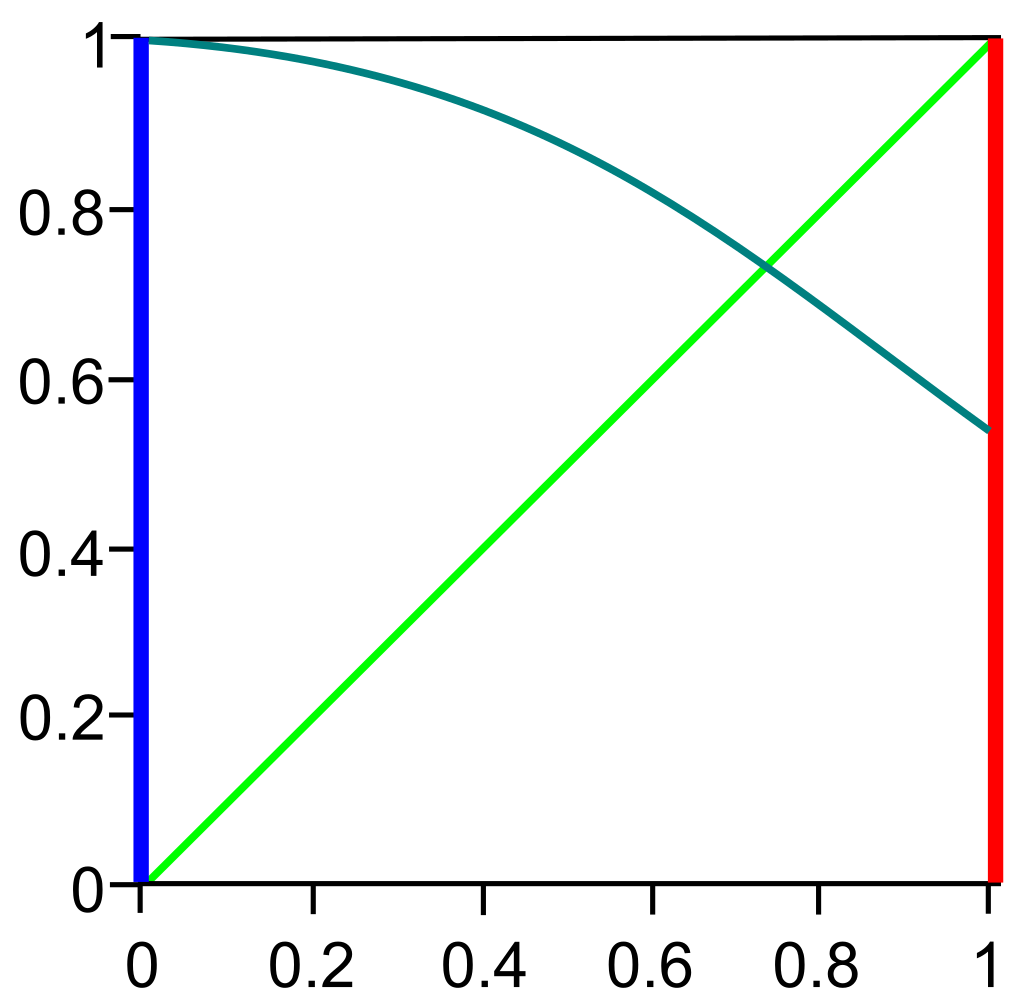}
    \caption{Intersection Between Light Green Line and Dark Green Curve}
    \label{fig:1d_brouwer}
\end{figure}

\textbf{Definition 1}: $B$ is homeomorphic to $A$ if there exists a continous bijecetion $f : A \rightarrow B$ such that $f^{-1}$ is also continuous, which is called a homeomorphism.

\textbf{Corollary 1.1}: For all $n \in \mathbb{N}$, if $K$ is homeomorphic to $\Delta_n$ and $f: K \rightarrow K$ is continuous, then $f$ has a fixed point.

\textbf{Theorem 1}: Every compact, convex subset of $R^m$ with nonempty interior is homeomorphic to a closed unit ball in Rm.

\textbf{Theorem 2}: Brouwer Fixed Point Theorem. Below are 2 equivalent (and True) definitions of BFP Theorem: (Refer to Appendix for proof of Brouwer Fixed Point Theorem using Sperner's Lemma)

(Br1) Every continuous map f: $B^d \rightarrow B^d$ has a fixed point. \cite{Ziegler}

(Br2) Every continuous map f: $B^d \rightarrow S^{d-1}$ has a fixed point. \cite{Ziegler}

\begin{proof} $(Br1) \iff (Br2)$:
$(Br1) \implies (Br2)$ This is trivial because $S^{d-1} \subseteq B^d$. For the converse, $(Br2) \implies (Br1)$: If $f: B^d \rightarrow B^d$ has no fixed point, then we set $g(x) \coloneq \frac{f(x)-x}{|f(x)-x|}$. This defines a map $g: B^d \rightarrow S^{d-1}$ that has a fixed point $x_{0} \in S^{d-1}$ by (Br2), with $x_{0} = \frac{f(x)-x}{|f(x)-x|}$. But this implies $f(x_0)=x_0(1+t)$ for $t \coloneq |f(x_0)-x_0|$, and this is impossible for $x_0 \in S^{d-1}$.
\end{proof}

\textbf{Theorem 2.1}: Brouwer Fixed Point Theorem: If $K$ os a nonempty, compact, convex subset of $\mathbb{R}^m$ and $f: K \rightarrow K$ is continuous, then $f$ has a fixed point.

\begin{proof}
We prove the result using strong induction on the dimension \( m \).

\textbf{Base Case: \( m = 0 \)}

In this case, the set \( K \) is a single point \( x^* \). Since \( f: K \rightarrow K \), it follows directly that
\[
f(x^*) = x^*,
\]
so \( x^* \) is a fixed point.

\textbf{Inductive Step}

Assume that the statement holds for all dimensions \( n < m \), i.e., for \( n = 0, 1, \ldots, m-1 \).

Now consider a compact, convex set \( K \subset \mathbb{R}^m \):

\begin{itemize}
    \item If \( K \) lies in a lower-dimensional affine subspace (i.e., it has at most \( n < m \) affinely independent vectors), then we can embed \( K \) in a lower-dimensional space where the inductive hypothesis applies. Therefore, \( f \) must have a fixed point.
    
    \item If \( K \) has \( m \) affinely independent vectors, then \( K \) has nonempty interior in \( \mathbb{R}^m \) due to convexity. By Theorem 1, \( K \) is homeomorphic to the closed unit ball \( \overline{B}(0,1) \subset \mathbb{R}^m \). 
    
    Thus, there exists a homeomorphism \( g: K \rightarrow \overline{B}(0,1) \), i.e., a continuous bijection with a continuous inverse.
    
    The standard simplex \( \Delta_{m-1} \) is also a compact, convex subset of \( \mathbb{R}^m \) with nonempty interior, so it too is homeomorphic to \( \overline{B}(0,1) \). Therefore, there exists a homeomorphism \( h: \overline{B}(0,1) \rightarrow \Delta_{m-1} \).
    
    Composing these maps yields a homeomorphism
    \[
    h \circ g: K \rightarrow \Delta_{m-1},
    \]
    with continuous inverse given by \( (h \circ g)^{-1} = g^{-1} \circ h^{-1} \).
    
    Consequently, \( K \) is homeomorphic to \( \Delta_{m-1} \). By Corollary 1.1, it follows that \( f \) has a fixed point.
\end{itemize}
\end{proof}

\section{Kakutani Fixed Point Theorem}
Kakutani’s Fixed Point Theorem builds on Brouwer’s Fixed Point Theorem by applying it to set-valued (multi-valued) functions. Like Brouwer’s theorem, it requires that the set $X$ be nonempty, compact, and convex. Under these conditions, if a set-valued function is upper semi-continuous and maps $X$ into closed, convex subsets of itself, then a fixed point is guaranteed to exist. Understanding such functions requires new definitions of continuity and fixed points, which are beyond the scope of this paper.

\textbf{Theorem 3}: Kakutani Fixed Point Theorem: If $S$ is a nonempty, compact, convex set in a Euclidean space and $\Phi: S \rightarrow \mathcal{P}(S)$ is upper semi-continuous, then $\Phi$ has a fixed point.

\section{Hex Theorem}

To prove the Hex Theorem, we begin with a fundamental result from \cite{MIT, Mark} in graph theory:

\textbf{Lemma 4.1.} Any finite graph where each vertex has degree at most two can be decomposed into a collection of disjoint subgraphs, each of which is either:  
(i) an isolated vertex,  
(ii) a simple cycle, or  
(iii) a simple path.

\begin{proof}
We proceed by induction on the number of edges in the graph.

Let $G$ be a finite graph with $N$ vertices, where each vertex has degree at most 2. Then, the maximum number of edges in $G$ is $N$. Denote a graph with $k$ edges as $G_k$.

\textbf{Base case:} When $k = 0$, the graph $G_0$ has no edges, and all vertices are isolated. Hence, the decomposition trivially holds.

\textbf{Inductive step:} Assume the lemma holds for any graph with $n$ edges, i.e., $G_n$ is a union of disjoint subgraphs each of which is an isolated vertex, simple cycle, or simple path.

Consider a graph $G_{n+1}$ with $n+1$ edges. Select and remove an arbitrary edge $(u,v)$. The resulting graph is $G_n$, which by the inductive hypothesis satisfies the lemma.

After removing the edge $(u,v)$, both vertices $u$ and $v$ have degree at most 1 (since their degree was at most 2 before). Therefore, $u$ and $v$ cannot lie on any cycle in $G_n$.

Now, reinsert the edge $(u,v)$ into $G_n$ to reconstruct $G_{n+1}$. The only subgraphs affected are those containing $u$ or $v$:
- If $u$ and $v$ were in separate paths, adding $(u,v)$ merges them into a longer path or creates a cycle.
- If one or both were isolated, adding $(u,v)$ forms a path.
- In all cases, the result is still a collection of disjoint isolated vertices, simple paths, or simple cycles.

Hence, the lemma holds for $G_{n+1}$.

By induction, the result holds for all $G_k$ with $0 \leq k \leq N$.
\end{proof}

 \begin{figure} [H]
    \centering
    \includegraphics[width=\linewidth]{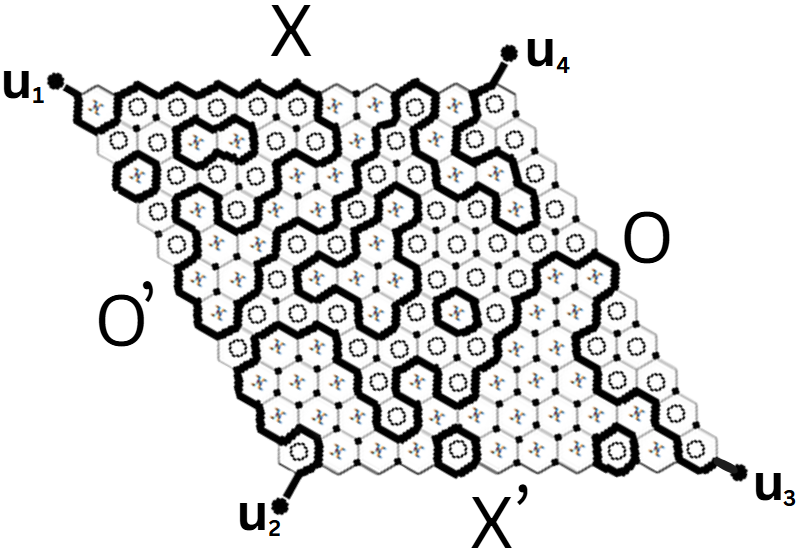}
    \caption{Hex board with the X and O notation where Red (Top and Bottom) and Blue (Left and Right)}
    \label{fig:hex_XO_notation}
\end{figure}

 For simplicity, we substitute the Red and Blue colors in the game with $x's$ and $o's$, respectively. We represent the game board as a graph $G = (V,E)$, with a set of nodes V and a set of edges E. Each vertex of a hexagonal board space is a node in V, and each side of a hexagonal board space is an edge in E. We create four additional nodes, one connected to each of the four corners of the core graph; call these new nodes $u_1, u_2, u_3,$ and $u_4$ and the edges that connect them to the core graph $e_1, e_2, e_3$ and $e_4$. An X-face is either a tile marked with an $x$ or one of the regions marked $X$ or $X^{\prime}$. Similarly, O-face is either a tile marked with an $o$ or one of the regions marked $O$ or $O^{\prime}$. Hence, the edges $e1, e2, e3, e4$ lie between an X-face and an O-face since the regions $X$, $X^{\prime}$, $O$, and $O^{\prime}$ are considered faces as well. Figure 4 shows a Hex board with the X and O notation. \cite{MIT}

 To prove the Hex Theorem, it suffices to show that two vertices out of $u1, u2, u3,$ and $u4$ are connected by a simple path. The hexagonal tiles traced out by this simple path contain a winning chain.

\textbf{Theorem 4 (Hex Theorem).}  
If every tile on a Hex board is marked either with $x$ or $o$, then either there exists an $x$-path connecting sides $X$ and $X'$ or an $o$-path connecting sides $O$ and $O'$.

\begin{proof}
We begin by constructing a subgraph $G' = (V, E')$ from the original Hex graph $G = (V, E)$, retaining the same set of vertices $V$ but selecting only certain edges to include in $E'$. An edge is included in $E'$ if it connects two adjacent hexagons marked with different symbols — one with $x$ and the other with $o$. Thus, only those edges that lie on the interface between $x$- and $o$-marked tiles are included.

As a result, the boundary vertices $u_1, u_2, u_3,$ and $u_4$ each have degree one in $G'$, since they lie between regions of different markings. Nodes entirely surrounded by tiles of the same mark will have degree zero. Nodes adjacent to two tiles of one mark and one tile of the opposite mark will have degree two. Therefore, all nodes in $G'$ have degrees of either 0, 1, or 2.

By a known graph-theoretical result, any graph where all vertices have degrees at most two consists of a collection of isolated vertices, simple cycles, and simple paths. Given that the boundary vertices $u_1$ to $u_4$ each have degree one, they must be endpoints of some simple paths.

Since the components of $G'$ are disjoint, the paths starting at $u_1$ through $u_4$ cannot form cycles. Thus, there must exist two simple paths in $G'$ connecting pairs among $u_1, u_2, u_3,$ and $u_4$. These paths form the interface separating regions marked with $x$ from those marked with $o$ and delineate a winning chain for one of the players.

For example, in Figure 4, \cite{MIT} the highlighted edges and vertices belong to $G'$, showing two disjoint simple paths — one connecting $u_1$ to $u_4$, and another connecting $u_2$ to $u_3$. These paths define a winning chain for the $o$-player.

Hence, regardless of the configuration of $x$ and $o$ on the board, there must always be a winning path for one of the two players.
\end{proof}

\begin{proof}[Proof (Strategy-Stealing Argument)]
\mbox{} Using the Hex Theorem, John Nash showed that the first player in Hex always has a winning strategy. The argument proceeds by contradiction and relies on a strategy-stealing approach: assume that the second player has a winning strategy. Then the first player makes an arbitrary initial move and pretends to be the second player from that point onward, following the supposed winning strategy. If the strategy ever requires playing on the already occupied cell, the first player simply makes another arbitrary move, which is always possible on the first move. This contradiction implies that the second player cannot have a winning strategy, and hence the first player must. Assume that the second player has a guaranteed winning strategy. Then, the first player can make an arbitrary move anywhere on the board and proceed to follow the second player's supposed winning strategy as if they were the second player. If the winning strategy ever instructs a move on the already occupied cell, the first player can simply choose another arbitrary cell for the first move — since having an extra piece on the board can never be a disadvantage. If the extra move is part of the winning path, it helps; if not, it does no harm. Therefore, the assumption that the second player has a winning strategy leads to a contradiction. This means that the first player must have a winning strategy.
\end{proof}

\section{Equivalence of the Hex Theorem and the Brouwer Fixed Point Theorem}

Before we begin, we adopt a new representation for the Hex board similar to \cite{MIT, Mark}.

Let $\mathbb{Z}^n$ represent the integer lattice points in $\mathbb{R}^n$. For distinct points $x, y \in \mathbb{R}^n$, define $|x - y| = \max_i |x_i - y_i|$, and say $x < y$ if $x_i \leq y_i$ for all $i$. Points $x$ and $y$ are comparable if $x < y$ or $y < x$.

Define the 2D Hex board of size $k$, denoted $B_k$, as a graph with vertices $z \in \mathbb{Z}^2$ such that $(1,1) \leq z \leq (k,k)$. Two vertices $z$ and $z'$ are adjacent if $|z - z'| = 1$ and they are comparable. We label the boundary edges as North (N), South (S), East (E), and West (W), where a vertex $z = (z_1, z_2)$ lies on:

- N if $z_2 = k$

- S if $z_2 = 0$

- E if $z_1 = k$

- W if $z_1 = 0$

We now replace the traditional players X and O with "horizontal" and "vertical" players.

We state the two theorems we wish to show are equivalent:

\textbf{Theorem 5 (Hex Theorem, Restated).} If $B_k$ is divided into two sets $H$ (horizontal) and $V$ (vertical), then either $H$ contains a connected path joining the E and W edges, or $V$ contains a connected path joining the N and S edges.

\textbf{Theorem 6 (Brouwer Fixed Point Theorem, Restated).} Any continuous function $f: I^2 \rightarrow I^2$ (where $I^2$ is the unit square) has a point $x$ for which $f(x) = x$.

\begin{proof}[Hex Theorem $\Rightarrow$ Brouwer Fixed Point Theorem]
Let $f: I^2 \rightarrow I^2$ be continuous, and suppose we want to find $x$ such that $|f(x) - x| < \epsilon$ for arbitrary $\epsilon > 0$. Since $f$ is uniformly continuous (as $I^2$ is compact), there exists $\delta > 0$ such that $|x - x'| < \delta$ implies $|f(x) - f(x')| < \epsilon$. Choose $k$ so large that $1/k < \delta$.

Define the following subsets of the $k \times k$ board $B_k$:
\begin{align*}
H^+ &= \{z \in B_k \mid f_1(z/k) - z_1/k > \epsilon\} \\
H^- &= \{z \in B_k \mid z_1/k - f_1(z/k) > \epsilon\} \\
V^+ &= \{z \in B_k \mid f_2(z/k) - z_2/k > \epsilon\} \\
V^- &= \{z \in B_k \mid z_2/k - f_2(z/k) > \epsilon\}
\end{align*}

We aim to show that not all points in $B_k$ lie in $H^+ \cup H^- \cup V^+ \cup V^-$. Any such point not covered would satisfy $|f(z/k) - z/k| < \epsilon$.

Suppose $z \in H^+$ and $z' \in H^-$ are adjacent. Then:
\begin{align*}
f_1(z/k) - z_1/k &> \epsilon \\
z_1'/k - f_1(z'/k) &> \epsilon
\end{align*}
Adding yields:
\[
f_1(z/k) - f_1(z'/k) + z_1'/k - z_1/k > 2\epsilon
\]
However, $|z_1/k - z_1'/k| \leq 1/k < \delta < \epsilon$, so this implies:
\[
f_1(z/k) - f_1(z'/k) > \epsilon
\]
Which contradicts the uniform continuity of $f$, so $H^+$ and $H^-$ are not contiguous. Similar logic applies to $V^+$ and $V^-$.

Let $H = H^+ \cup H^-$ and $V = V^+ \cup V^-$. Any connected set in $H$ must be entirely in either $H^+$ or $H^-$, since the two are not contiguous. But $H^+$ cannot reach the E boundary (since $f_1(z/k) \leq 1$), and $H^-$ cannot reach the W boundary. Thus, no path in $H$ connects E to W. Similar reasoning shows $V$ cannot connect N to S. 

Therefore, by the Hex Theorem, $H \cup V$ cannot cover all of $B_k$, and any point not in $H \cup V$ must satisfy $|f(z/k) - z/k| < \epsilon$.
\end{proof}

\begin{proof}[Brouwer Fixed Point Theorem $\Rightarrow$ Hex Theorem]
We assume $B_k$ is covered by two disjoint sets $H$ and $V$ such that:

- No path in $H$ connects E to W

- No path in $V$ connects N to S

Let $\widehat{W}$ be all vertices in $H$ connected to the W boundary, and $\widehat{E} = H \setminus \widehat{W}$. Similarly, let $\widehat{S}$ be vertices in $V$ connected to S, and $\widehat{N} = V \setminus \widehat{S}$.

Define $f: B_k \rightarrow B_k$ by:
\[
f(z) =
\begin{cases}
z + e^1, & z \in \widehat{W} \\
z - e^1, & z \in \widehat{E} \\
z + e^2, & z \in \widehat{S} \\
z - e^2, & z \in \widehat{N}
\end{cases}
\]
This function is well-defined, since no set in the domain borders its complement in a way that would violate board bounds.

Extend $f$ to a continuous simplicial map $\widehat{f}$ on $I_k^2$ using convex combinations of triangle vertices. Each triangle's image under $\widehat{f}$ is determined by vector displacements within a single quadrant. By a lemma: such a map has no fixed point if the origin is not in the convex hull of displacement vectors.

Thus, $\widehat{f}$ has no fixed point, which contradicts Brouwer's Theorem. Therefore, the assumption must be false, and either $H$ connects E to W or $V$ connects N to S.
\end{proof}

\section{Conclusion}
Brouwer's Fixed Point Theorem, Kakutani's Fixed Point Theorem and its connection with Combinatorics and Geometry seem not to be as well known as the applications of the Borsuk-Ulam Theorem. In this paper, we delve into Brouwer Fixed Point Theorem and its connection with Hex Board Game. We also see how Hex Theorem is equivalent to Brouwer Fixed Point Theorem to offer a more intuitive approach to understanding the Fixed Point Theorems through Board Games like Hex. Fixed Point Theorems are powerful and key concepts used to prove the existence of Nash Equilibria in games.

\section*{Acknowledgments}

I would like to express my gratitude to my supervisor, Dr Fedor Duzhin, I am grateful to be able to work on a project of my interest and under the expertise and supervision of Dr Fedor Duzhin.

\textbf{All papers arising as a result of URECA programme must have the following (important) statement under the acknowledgement section.}

\textit{I would like to acknowledge the funding support from Nanyang Technological University – URECA Undergraduate Research Programme for this research project.}

\textbf{Acknowledgement also applies to Journal papers, Conference papers and Proceedings of URECA Undergraduate Research papers arising as a result of URECA.}

%\section{About references}
%Then cite as follows: \cite{AllenYuan}. To
%cite several references together, format the citation as %\cite{Ziegler, Mark, MIT}.

% \bibliographystyle{plain}  % Use an appropriate style (plainnat, abbrvnat, etc.)
% \bibliography{references}  % This refers to references.bib (without .bib extension)

\section*{Appendix}
We begin with the following definitions, building toward Sperner’s Lemma, which plays a central role in proving the Brouwer Fixed Point Theorem \cite{AllenYuan}.

\textbf{Definition A1.1.} A subset $S \subseteq \mathbb{R}^m$ is \textit{convex} if for all $x, y \in S$ and all $\lambda \in [0, 1]$, we have $\lambda x + (1 - \lambda)y \in S$.

\textbf{Definition A1.2.} An \textit{$n$-simplex} is the set of all strictly positive convex combinations of $n+1$ affinely independent points. Specifically, for vertices $x^0, \ldots, x^n$, we define:
\[
T = \left\{ \sum_{i=0}^{n} \lambda_i x^i \ \middle| \ \lambda_i > 0 \ \forall i, \ \sum_{i=0}^{n} \lambda_i = 1 \right\}
\]

The \textit{standard $n$-simplex} is the one formed by the $n+1$ standard basis vectors in $\mathbb{R}^{n+1}$.

\textbf{Definition A1.3.} The \textit{closed $n$-simplex} is the convex hull of the vertices $x^0, \ldots, x^n$:
\[
\overline{T} = \left\{ \sum_{i=0}^{n} \lambda_i x^i \ \middle| \ \lambda_i \geq 0, \ \sum_{i=0}^{n} \lambda_i = 1 \right\}
\]

\textbf{Definition A1.4.} A \textit{simplicial subdivision} of $\overline{T}$ is a finite collection of simplices $\{T_i\}_{i=1}^{m}$ such that:
\begin{multline*}
\bigcup_{i=1}^{m} T_i = \overline{T}, \quad \text{and for all } i, j, \ \overline{T}_i \cap \overline{T}_j \\
\text{is either empty or a common face of both.}
\end{multline*}

\textbf{Definition A1.5.} Let $y = \sum_{i=0}^{n} \lambda_i x^i$ lie in the convex hull of $\{x^0, \ldots, x^n\}$. Define:
\[
\chi(y) = \{i \in \{0, \ldots, n\} \mid \lambda_i > 0\}
\]
Then $y$ belongs to the simplex spanned by $\{x^i \mid i \in \chi(y)\}$.

\textbf{Definition A1.6.} Let $V$ be the set of all vertices of the simplicial subdivision of $\overline{T}$. A labeling function $f: V \to \{0, \ldots, n\}$ is a \textit{proper labeling} if for every $v \in V$, we have $f(v) \in \chi(v)$. A simplex is \textit{completely labeled} if its vertices receive all labels $0, \ldots, n$.

\textbf{Theorem A1 (Sperner’s Lemma).} Let $\overline{T}$ be a closed $n$-simplex that is simplicially subdivided and properly labeled. Then there exists an odd number of completely labeled $n$-subsimplices — in particular, at least one.

\textbf{Lemma A1.7.} If $x, y \in \Delta_m$ and $x \leq y$ component-wise, then $x = y$.

\section*{Theorem A1.1 (Brouwer Fixed Point Theorem)}
Any continuous function $f: \Delta_m \rightarrow \Delta_m$ has a fixed point.

\begin{proof}
Let $f: \Delta_m \to \Delta_m$ be continuous. For each $n \in \mathbb{N}$, take a simplicial subdivision of $\Delta_m$ such that the \textit{mesh} — the maximum diameter of any subsimplex — is at most $1/n$. Such subdivisions exist, e.g., via barycentric subdivision.

Let $V$ be the set of all vertices in the subdivision. Define a labeling function $\lambda: V \to \{0, \ldots, m\}$ as follows. For each $v \in V$, define:
\[
\lambda(v) = \min \left\{ i \in \chi(v) \mid f_i(v) \leq v_i \right\}
\]
This minimum exists: if no such $i$ exists, then
\[
\sum_{i=0}^{m} f_i(v) > \sum_{i=0}^{m} v_i = 1,
\]
contradicting that $f(v) \in \Delta_m$. Thus, $\lambda$ is a proper labeling.

By Sperner's Lemma, there exists a completely labeled $m$-simplex in the subdivision. Let its vertices be ${}^n p^0, \ldots, {}^n p^m$. Repeat this process for each $n$, producing a sequence of such simplices.

Because the mesh goes to zero, the diameter of these simplices also tends to zero. So the vertices ${}^n p^0, \ldots, {}^n p^m$ get arbitrarily close. Since $\Delta_m$ is compact, we can extract a convergent subsequence from each vertex sequence. Using a diagonal argument, we obtain a point $p \in \Delta_m$ such that:
\[
\lim_{n \to \infty} {}^n p^i = p \quad \text{for all } i.
\]
From the labeling rule, we have $f_i({}^n p^i) \leq {}^n p^i$ for all $n$ and all $i$. By continuity of $f$, this implies:
\[
f_i(p) \leq p_i \quad \forall i.
\]
So $f(p) \leq p$, and since both $f(p)$ and $p$ lie in $\Delta_m$, Lemma A1.7 yields $f(p) = p$. Hence, $p$ is a fixed point.
\end{proof}

\section*{Remark}
The Brouwer Fixed Point Theorem holds more generally: any continuous function $f: X \to X$ on a nonempty, compact, convex subset $X \subseteq \mathbb{R}^m$ has a fixed point. This follows because such a set $X$ is homeomorphic to the standard closed $m$-simplex, and homeomorphisms preserve the fixed-point property.

\end{multicols}

\end{document}